\documentclass[11pt]{article}
\usepackage{amsmath, amsthm, amssymb}    % May not all be necessary
\usepackage{bridges}			% Custom bridges proceedings style
\usepackage{graphicx}			% For including pictures
\usepackage[colorlinks=true, urlcolor=blue, citecolor=black, linkcolor=black]{hyperref}  
\usepackage{subcaption}			% For captioning multi-panel figures

\usepackage{float}
\usepackage{enumitem}
\newcommand{\R}{\mathbb{R}}
\renewcommand{\phi}{\varphi}
\newcommand{\euc}{\mathcal{E}}
\usepackage{tikz}

\title{The Flat Klein Bottle Rendered in Curved-Crease Origami}

\author{Stepan Paul
\vspace{10pt}\\
North Carolina State University; sspaul2@ncsu.edu\\
} % end \author

\date{}					% Suppress any date on submissions

\begin{document}

\maketitle

% Prevent page number 1 from being printed on the first page.
\thispagestyle{empty}

\begin{abstract}

We introduce a simple and concrete way of visualizing in three dimensions a ``flat'' Klein bottle---one whose local intrinsic geometry is the same as that of a flat plane---which preserves most its topological and geometric structure. Concretely, the flatness property means that a small patch of the surface around any point can be flattened to a patch of the plane without stretching or compressing. Thus we can use the medium of curved-crease origami with inelastic film to make a model which, except for its self-intersections, necessarily has the flatness property, even along its folded edges. As such, the sculpture presented here illustrates both the flatness, and, through its coloring, the non-orientability of a Klein bottle.

\end{abstract}

% Bridges papers are usually no more than 8 pages in length.  So
% there's really no need to have numbered sections, unless the
% author really needs to refer to sections by number within the paper's text.  
% So to suppress sequential section numbers, append an asterisk to 
% the \section command, as in:

%%%%%%%%%%%%%%%%%%%%%%%%%%%%%%%%%%%%%%%%%%%
\section*{Introduction}

We present here a sculpture of a flat Klein bottle meant to illustrates both its topological and geometric properties.

\begin{figure}[H]
 \begin{center}
  \includegraphics[width=4in]{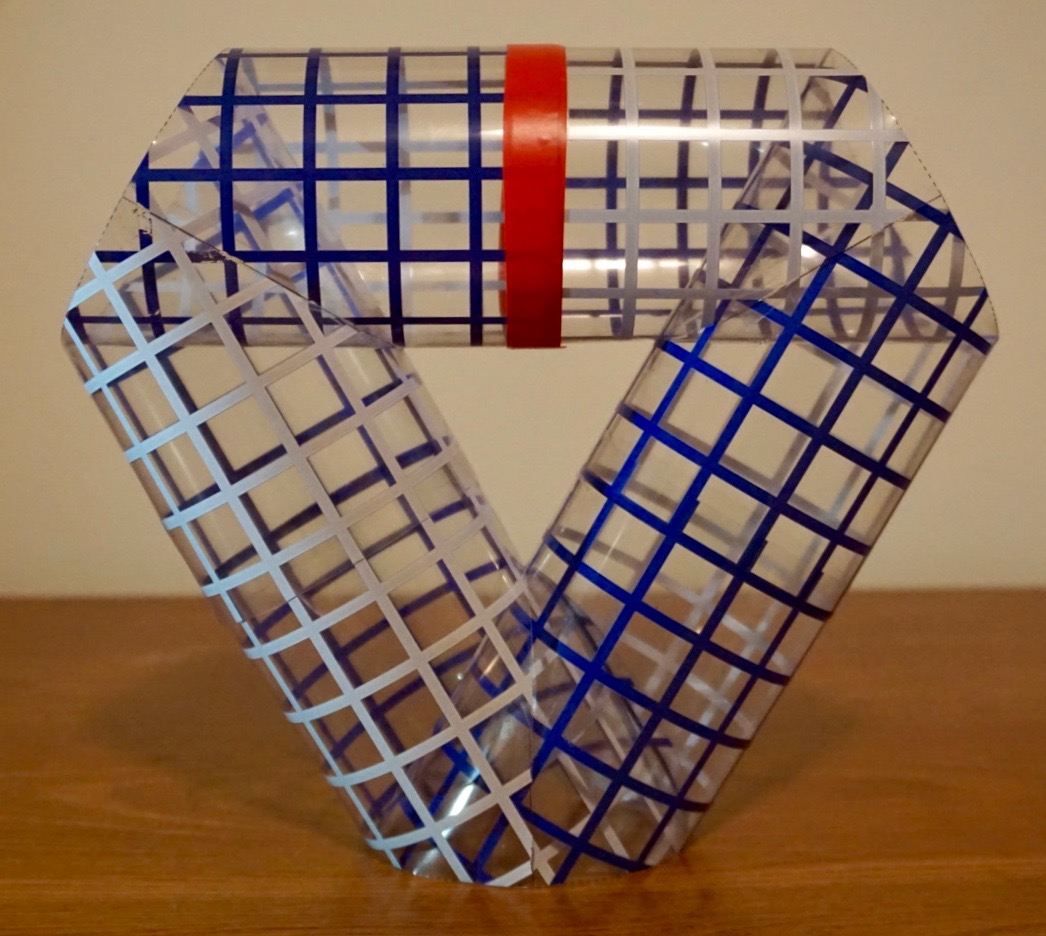}
  \caption{A curved-crease origami sculpture of a flat Klein bottle made from printed transparency film.}
 \end{center}
\end{figure}

Here the term ``flat'' is used to mean that the surface has the same local intrinsic geometry as a flat plane at every point. I elaborate on this property and its significance in the next section, but intuitively, it means that a patch of the surface around any point can be deformed into a flat plane without stretching or compressing---that is, without distorting distances.

%While a flat surface need not \emph{be} a plane,  For a smooth surface, flatness means that the so-called ``Gaussian curvature'' is zero, but a non-smooth surface, like the one represented by the present sculpture can also be said to be flat if edges and vertices can be similarly deformed

In the present sculpture, the uniform grid pattern and origami medium are intended to illustrate the %existence of a Riemannian metric with everywhere-zero Gaussian curvature---which we will call a \emph{flat} metric---on the 
 flatness of Klein bottle, and the coloring of the grid is intended to illustrate its non-orientability.

We think of the underlying mathematical figure $\nabla$, which consists of three circular cylinders of equal radius centered around the axes of an equilateral triangle, represented by the sculpture as the image of the abstract flat Klein bottle under a map $\phi$ into euclidean three-space with the following properties.
\begin{enumerate}[label=\Roman*., ref=\Roman*]
 \item\label{property:c0} $\phi$ is continuous everywhere, and $\phi$ is smooth except along a finite set of curves,
 \item\label{property:local-isometry} $\phi$ is locally geometry-preserving (a 
 {local isometry}) except at a finite set of points, and
 \item\label{property:injective} $\phi$ is one-to-one, except along a finite set of curves.
\end{enumerate}

 It is not possible to embed a Klein bottle into euclidean three-space, and thus neither the topology nor the geometry can be completely preserved. However, this set of properties for $\phi$ can be thought of as an approximation to topology- and geometry-preservation. To the author's knowledge, this is the first explicit example of a map from the flat Klein bottle into euclidean three-space satisfying these conditions%\footnote{Although the existence of maps satisfying an even stronger set of conditions can be shown to exist using the results of Nash and Kuiper \cite{nash54,kuiper55}}
, and perhaps the first attempt to concretely visualize flat geometry on the Klein bottle in three dimensions.
 
 In what follows, I start with an overview of the mathematical concepts at play and prove that the claimed properties of $\phi$ follow from the theory of curved-crease origami. I then elaborate on the construction and intention of the sculpture itself, focusing on how it self-evidentially illustrates important properties of the Klein bottle. Finally, I then give some mathematical and artistic context by exploring the relationship of the present illustration to related works.
 
\section*{Mathematical Overview}

In this section, we explain in a bit more detail what is meant by ``flatness'', how it relates to the Klein bottle and torus specifically, and how the theory of curved-crease origami is applied in constructing the figure $\nabla$.

\subsection*{Flatness}

\begin{figure}[h!tbp]
\centering
\begin{minipage}[c]{0.3\textwidth} 
	\includegraphics[height=1.75in]{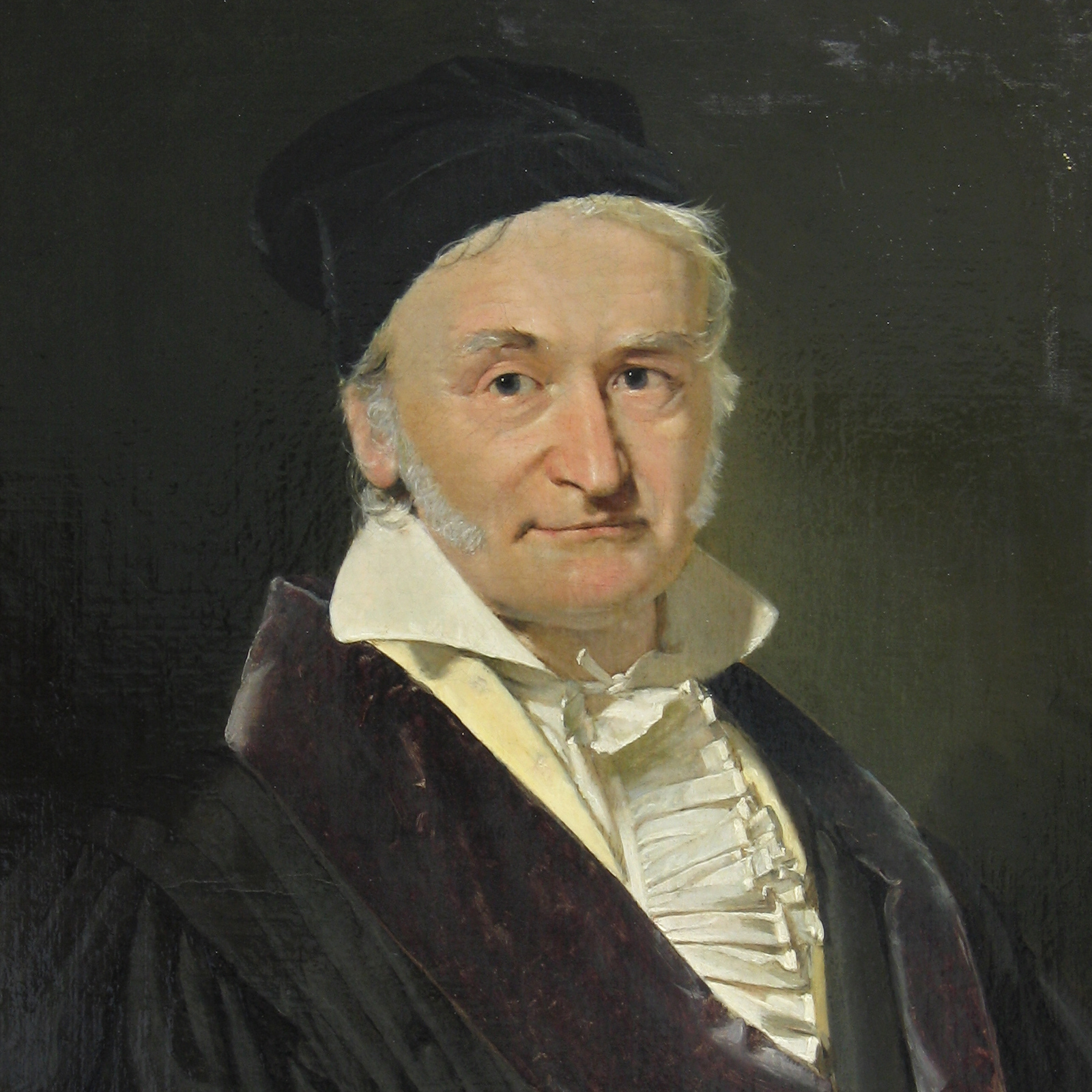}
        	\subcaption{} % Add subcaption text if desired, or use \subcaption* to suppress (a), (b), etc. labels
        	\label{fig:gauss}
\end{minipage}
~ %add desired spacing between images, e. g., ~, \quad, \qquad, \hfill etc.	
\begin{minipage}[c]{0.25\textwidth} 
	\includegraphics[height=1.75in]{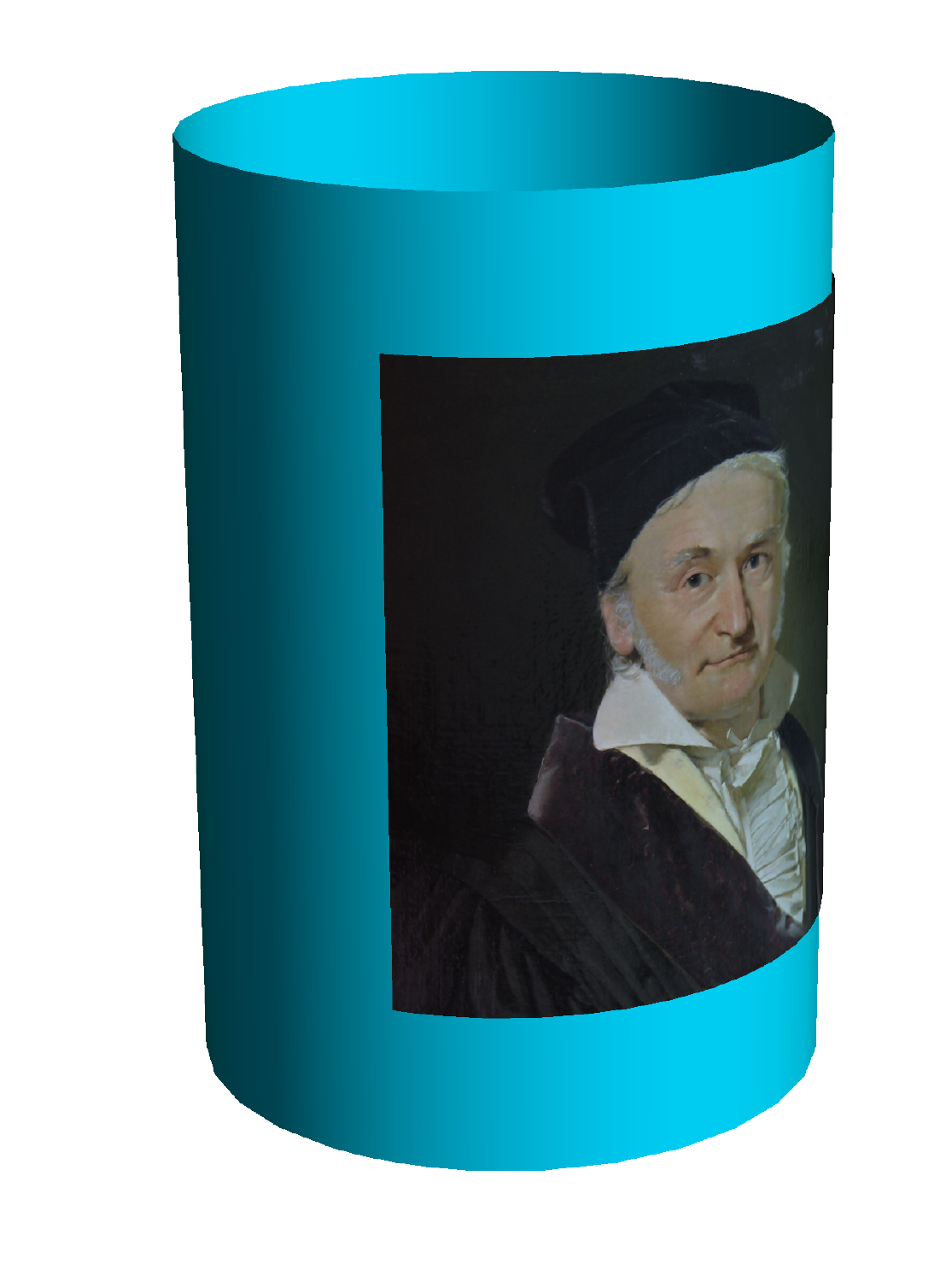}
        	\subcaption{} % Add subcaption text if desired, or use \subcaption* to suppress (a), (b), etc. labels
        	\label{fig:gauss-cylinder}
\end{minipage}
~ %add desired spacing between images, e. g., ~, \quad, \qquad, \hfill etc.	
\begin{minipage}[c]{0.35\textwidth} 
	\includegraphics[height=1.75in]{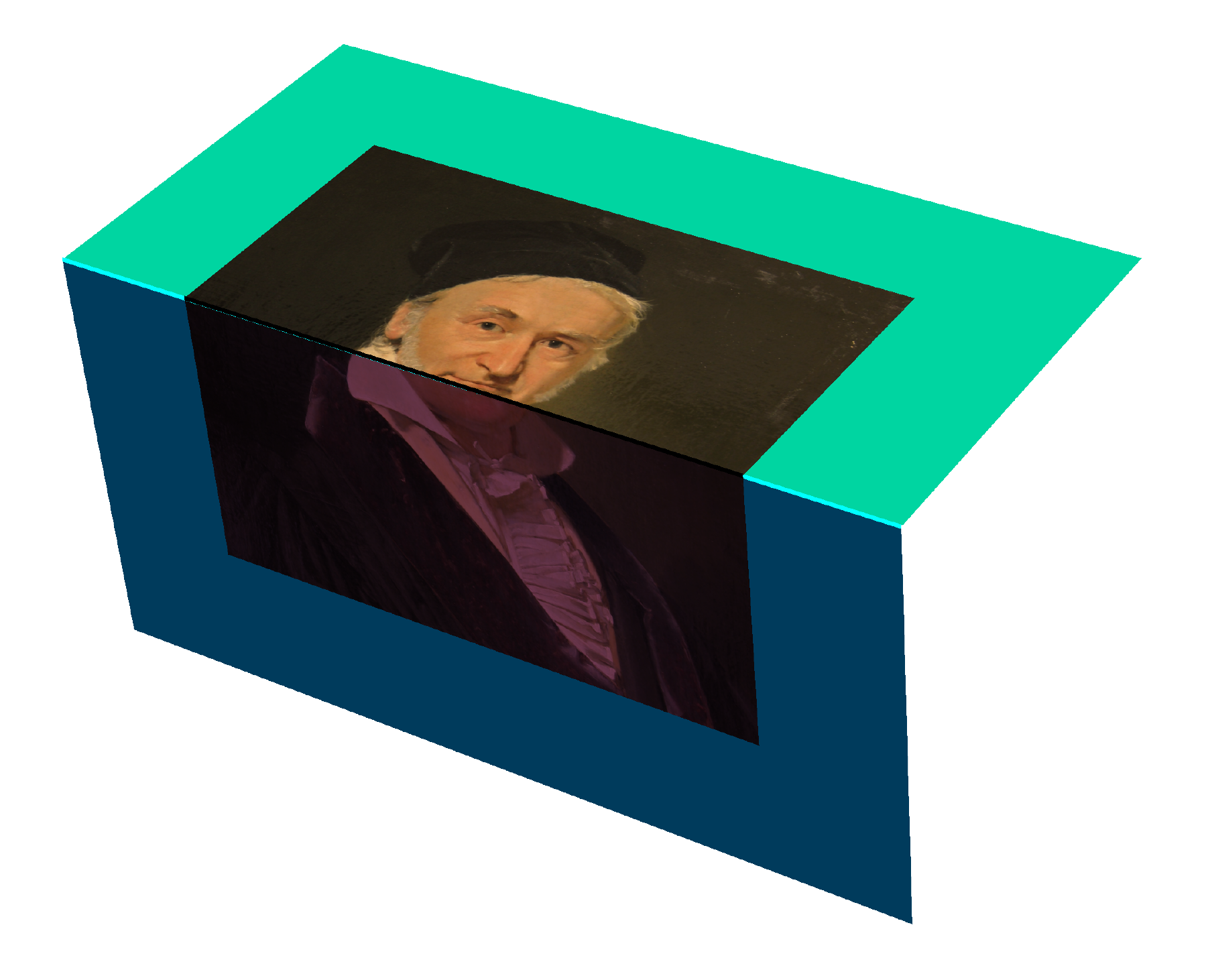}
        	\subcaption{} % Add subcaption text if desired, or use \subcaption* to suppress (a), (b), etc. labels
        	\label{fig:gauss-edge}
\end{minipage}
~ %add desired spacing between images, e. g., ~, \quad, \qquad, \hfill etc.	

\caption{Examples of flat surfaces: A piece of the euclidean plane (decorated with an image) maps isometrically (with the image undistorted) onto a cylinder (b) and a non-smooth surface (c).}
\label{flat-examples}
\end{figure}

As we mention in the introduction, a surface is said to be \emph{flat} if it has the same local intrinsic geometry as the euclidean plane. This means that for any point on the surface, a piece of the euclidean plane can map \emph{isometrically} onto a neighborhood of the point (i.e.~the map is such that the distance between two points in the plane is the same as the distance between their images \emph{along the surface}). Examples of flat surfaces include the cylinder and the non-smooth surface consisting of two half-planes meeting at any nonzero dihedral angle along a line, as illustrated in Figure \ref{flat-examples}.

By Gauss' \emph{Theorema Egregium}, a smooth flat surface can be characterized extrinsically by the fact that its Gaussian curvature is everywhere zero, but as in Figure \ref{fig:gauss-edge}, a flat surface need not be smooth.

In contrast, a spherical surface is not flat since it cannot be flattened without distorting distances. Likewise, the surface of a cube and the surface of a solid cylinder are not flat even though their faces are; neither a neighborhood of a vertex of the cube, nor a neighborhood of a point along the edge of the cylinder can flattened without stretching or ripping. 

\begin{figure}[h!tbp]
\centering
\begin{minipage}[c]{0.3\textwidth} 
	\includegraphics[height=1.75in]{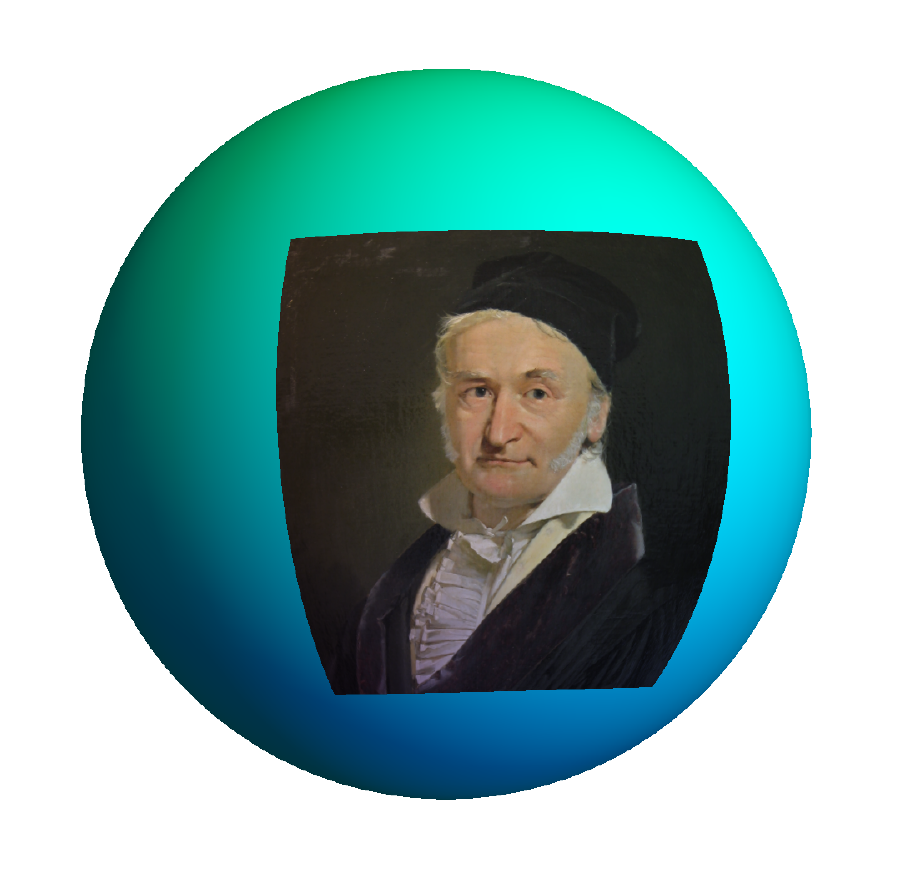}
        	\subcaption{} % Add subcaption text if desired, or use \subcaption* to suppress (a), (b), etc. labels
        	\label{fig:gauss-sphere}
\end{minipage}
~ %add desired spacing between images, e. g., ~, \quad, \qquad, \hfill etc.	
\begin{minipage}[c]{0.3\textwidth} 
	\includegraphics[height=1.75in]{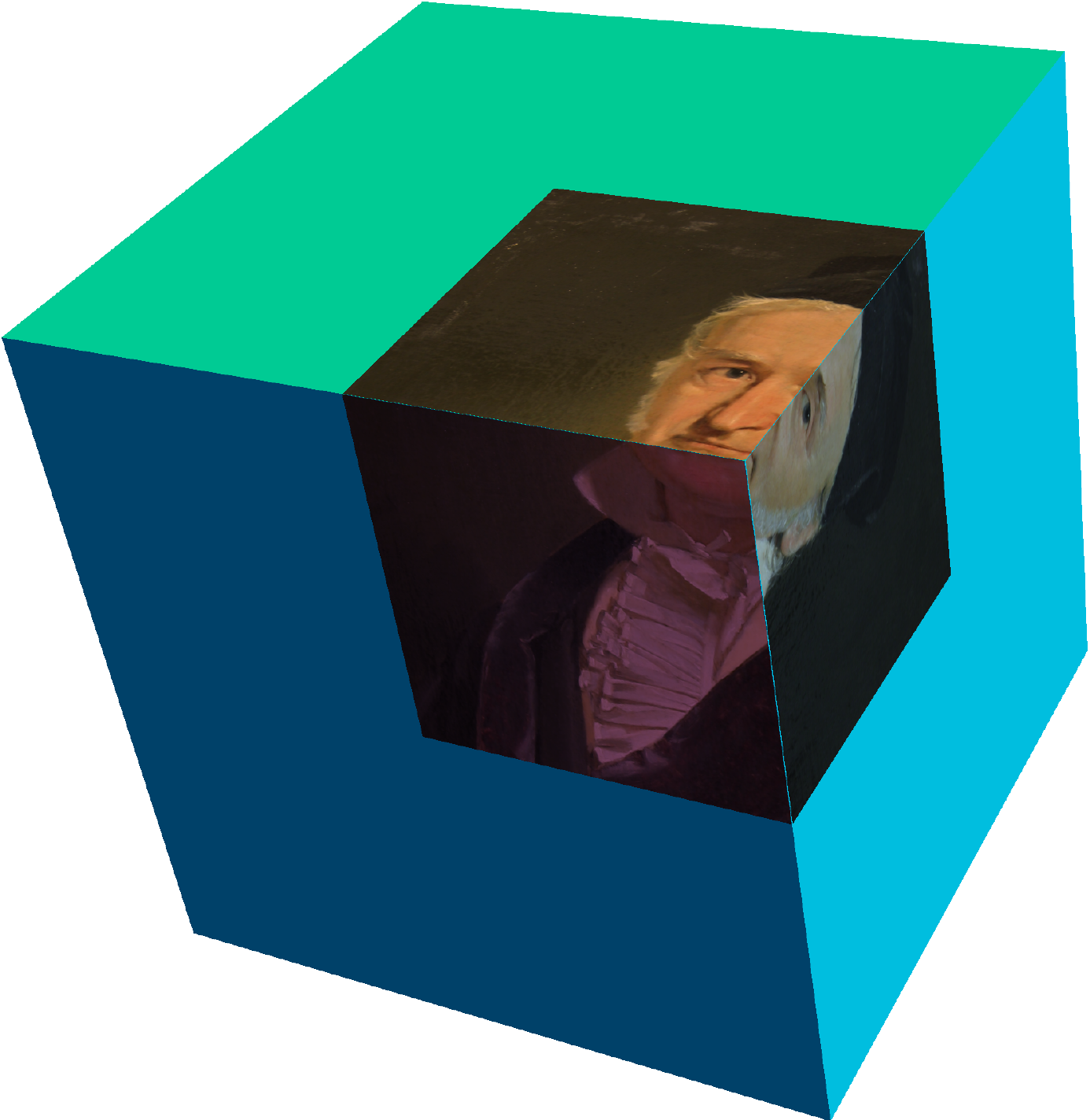}
        	\subcaption{} % Add subcaption text if desired, or use \subcaption* to suppress (a), (b), etc. labels
        	\label{fig:gauss-cube}
\end{minipage}
~ %add desired spacing between images, e. g., ~, \quad, \qquad, \hfill etc.	
\begin{minipage}[c]{0.35\textwidth} 
	\includegraphics[height=1.75in]{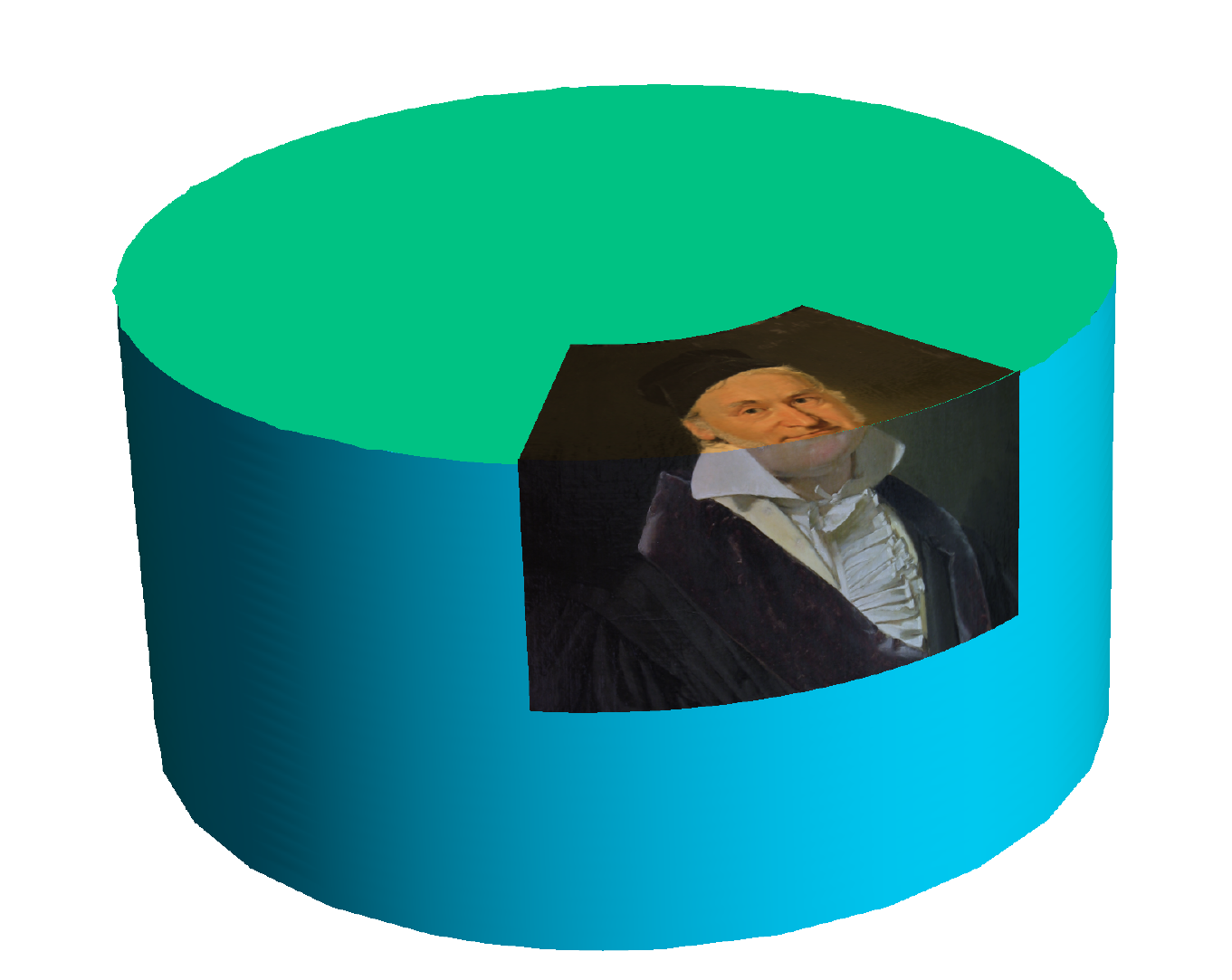}
        	\subcaption{} % Add subcaption text if desired, or use \subcaption* to suppress (a), (b), etc. labels
        	\label{fig:gauss-can}
\end{minipage}
~ %add desired spacing between images, e. g., ~, \quad, \qquad, \hfill etc.	

\caption{Examples of non-flat surfaces: A piece of the euclidean plane cannot be mapped (a) to a sphere, (b) onto a neighborhood of a vertex of a cube, nor (c) onto a neighborhood of a point along the edge of the surface of the solid cylinder without distorting distances.}
\label{non-flat-examples}
\end{figure}

In a way, the obstruction to flatness for all three of these examples is the same; their global topology---in each case, that of a sphere---restricts the geometry. Indeed, it is a consequence of the Gauss-Bonnet Theorem that any flat compact surface (without boundary) is topologically either a torus or Klein bottle.

As we discuss in the section on the mathematical and artistic context of this work, the flat torus has been illustrated in several previous works, and the present sculpture seeks to do the same for the flat Klein bottle.

\subsection*{The Flat Klein Bottle}

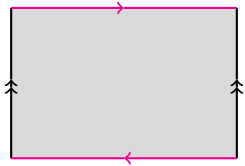
\begin{figure}[h!tbp]
 \begin{center}
 \begin{tikzpicture}
  \fill[black!15!white] (0,0) rectangle (3,2);
  
  \draw[thick,black,->] (0,0) -- (0,0.95);
  \draw[thick,black,->] (0,0.95) -- (0,1.05);
  \draw[thick,black] (0,1.05) -- (0,2);
  \draw[thick,black,->] (3,0) -- (3,0.95);
  \draw[thick,black,->] (3,0.95) -- (3,1.05);
  \draw[thick,black] (3,1.05) -- (3,2);
  
  \draw[thick,magenta,->] (3,0) -- (1.5,0);
  \draw[thick,magenta] (1.5,0) -- (0,0);
  \draw[thick,magenta,->] (0,2) -- (1.5,2);
  \draw[thick,magenta] (1.5,2) -- (3,2);
  
 \end{tikzpicture}
 \end{center}
 \caption{Edge identifications on a euclidean rectangle to form a flat Klein bottle.}
 \label{fig:klein-bottle-rectangle}
\end{figure}

A flat Klein bottle can be constructed by identifying the edges of a euclidean rectangle as show in Figure \ref{fig:klein-bottle-rectangle}.

A locally isometric map on the rectangle that respects the edge identifications can then be interpreted as a local isometry on the Klein bottle.

\subsection*{Theory of Curved-Crease Origami}

In constructing the map $\phi$ and its image $\nabla$, we apply the theory of curved-crease origami found, for example in \cite{ft99}. Denote by $\euc^2$ and $\euc^3$ the euclidean plane and three-space respectively. Let $\alpha$ be a curve in $\euc^2$---thought of as a curve drawn on a flat piece of paper---and let $f_1$ be a smooth local isometry of a neighborhood of $\alpha$ into $\euc^3$---thought of as a particular positioning of the paper in $\euc^3$, possibly bending (but not folding) the paper.

 The theory then says that if for every point $x$ along $\alpha$, the curvature of $f_1(\alpha)$ at $f_1(x)$ is strictly greater than the curvature of $\alpha$ at $x$, then there exists a unique smooth local isometry $f_2$ that is different from $f_1$ but which also maps a neighborhood of $\alpha$ into $\euc^3$ and which agrees with $f_1$ along $\alpha$. We interpret this to mean that there is exactly one other way to position the sheet of paper that places the drawn curve $\alpha$ in the same location that $f_1$ did.
 %(PICTURE: a) alpha drawn on a paper, b) f(alpha), g(alpha), c) piecewise)
 Furthermore, at points where the curvature of $\alpha$ is non-vanishing, the piecewise-defined function $f$ defined to be $f_1$ on one side of $\alpha$ and $f_2$ on the other side is also a local isometry.
 
 We apply this theory to the example of the elliptical intersection of a circular cylinder and a plane. In particular, if $\Psi_1$ is a circular cylinder of radius $r$ and $\Pi$ is a plane meeting the axis of $\Psi_1$ at an angle of $0<\theta<\frac{\pi}{2}$, then one can choose a coordinate system so that the map
\begin{equation}\label{eqn:cylinder-formula}
 f_1(u,v)=\left(r\sin\frac{u}{r},v,r\cos\frac{u}{r}\right)
\end{equation}
maps (locally isometrically) onto $\Psi$ and so that $\Pi$ has equations $y=\tau x$ where $\tau=\tan\theta$.

The pre-image of the the elliptical intersection $\epsilon$ is the curve $\alpha$ in $\R^2$ defined by $v=\tau r\sin\frac{u}{r}$. Applying the theorem locally, $f_1|\alpha$ extends to $f_1$ and exactly one other smooth local isometry $f_2$, and we can immediately identify $f_2$ as the reflection $\rho\circ f_1$, where $\rho$ is reflection across $\Pi$ (since $\rho$ fixes $\epsilon$), and thus the piecewise function
\begin{equation}\label{eqn:piecewise}
f(u,v)=\begin{cases} f_1(u,v) & v\leq \tau\sin\frac{u}{r}\\f_2(u,v) & v\geq \tau\sin\frac{u}{r}\end{cases}
\end{equation}
is a local isometry away from the inflection points of $\alpha$.

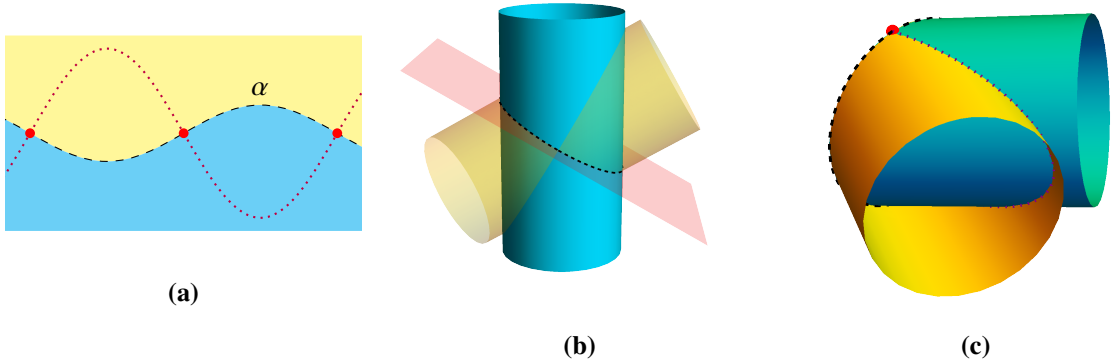
\begin{figure}[h!tbp]
\centering
\begin{minipage}[c]{0.3\textwidth} 
 \begin{center}
	\begin{tikzpicture}[scale=0.65]
	 \def\a{0.577}; % 1/sqrt(3)
   \def\b{1.732}; % sqrt(3)
   
   \clip (-pi-0.5,-2) rectangle (pi+0.5,2);

	 \fill[yellow!50!white] (-3*pi/2,\a) cos (-pi,0) sin (-pi/2,-\a) cos (0,0) sin (pi/2,\a) cos (pi,0) sin (3*pi/2,-\a) -- (3*pi/2,2) -- (-3*pi/2,2) -- (-3*pi/2,0);
	 
	 \fill[cyan!50!white] (-3*pi/2,\a) cos (-pi,0) sin (-pi/2,-\a) cos (0,0) sin (pi/2,\a) cos (pi,0) sin (3*pi/2,-\a) -- (3*pi/2,-2) -- (-3*pi/2,-2) -- (-3*pi/2,0);
	 
	 \draw[dashed] (-3*pi/2,\a) cos (-pi,0) sin (-pi/2,-\a) cos (0,0) sin (pi/2,\a) node[anchor=south]{$\alpha$} cos (pi,0) sin (3*pi/2,-\a);
	 
	 \draw[purple,thick,dotted] (-3*pi/2,-\b) cos (-pi,0) sin (-pi/2,\b) cos (0,0) sin (pi/2,-\b) cos (pi,0) sin (3*pi/2,\b);
	 
	 \fill[red] (-pi,0) circle (0.1cm);
	 \fill[red] (0,0) circle (0.1cm);
	 \fill[red] (pi,0) circle (0.1cm);

	\end{tikzpicture}
\end{center}
	
        	\subcaption{} % Add subcaption text if desired, or use \subcaption* to suppress (a), (b), etc. labels
        	\label{fig:joint-domain}
\end{minipage}
~ %add desired spacing between images, e. g., ~, \quad, \qquad, \hfill etc.	
\begin{minipage}[c]{0.3\textwidth} 
	\includegraphics[height=1.75in]{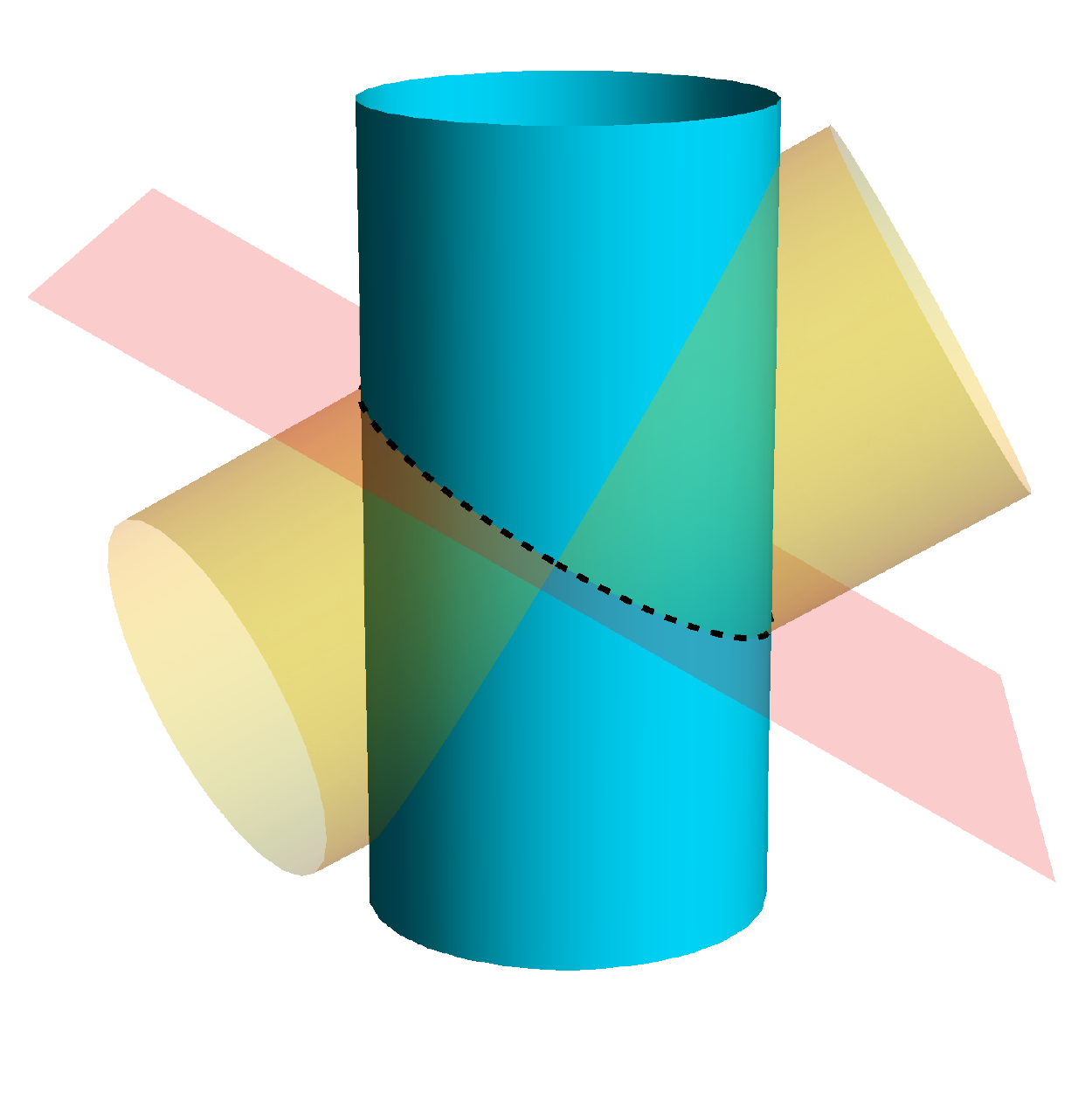}
        	\subcaption{} % Add subcaption text if desired, or use \subcaption* to suppress (a), (b), etc. labels
        	\label{fig:two-cylinders}
\end{minipage}
~ %add desired spacing between images, e. g., ~, \quad, \qquad, \hfill etc.	
\begin{minipage}[c]{0.3\textwidth} 
	\includegraphics[height=1.75in]{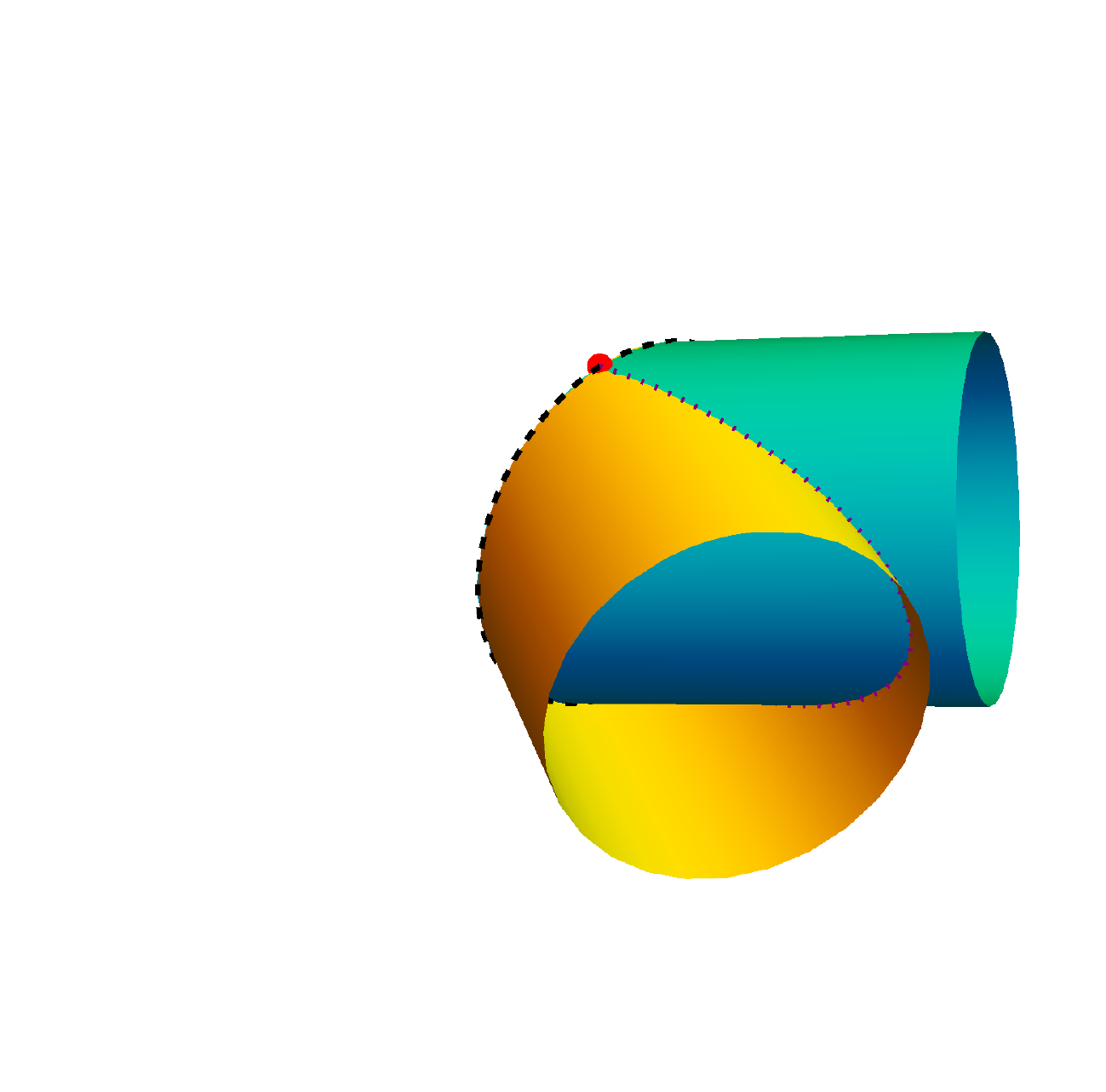}
        	\subcaption{} % Add subcaption text if desired, or use \subcaption* to suppress (a), (b), etc. labels
        	\label{fig:joint}
\end{minipage}
~ %add desired spacing between images, e. g., ~, \quad, \qquad, \hfill etc.	

\caption{The $uv$-plane in (a) is mapped by $f_1$ onto the solid cylinder in (b). The dashed curve $\alpha$ in (a) maps to the dashed intersection of the solid cylinder with a plane in (b). The image of $f_2$ is shown as a transparent cylinder in (b). The image of the piecewise function $f$ is shown in (c), and the inflection points of $\alpha$ and their images are highlighted, and the semi-ellipse where the two cylinders intersect and its pre-image are shown as dotted curves in (c) and (a) respectively.}
\label{Cylindrical joint}
\end{figure}

One further observation to make here is that the ``handedness'' of the cylinder parameterizations $f_1$ and $f_2$ are different in the following sense. If, while one points the thumb of one's right hand in the direction of increase in the axial variable ($v$) one's fingers curl in the direction of increase in the angular variable ($u$), we say that the cylinder parameterization is ``right-handed''; otherwise it is left-handed. Mathematically, the normal vector resulting from the cross product of the partial derivatives $\frac{\partial f}{\partial u}\times \frac{\partial f}{\partial v}$ points outward for a right-handed parameterization and inward for a left-handed parameterization; the reversal results from the fact that $f_2$ is a reflection of $f_1$.

\subsection*{Constructing $\nabla$ and $\phi$}

One can now imagine patching together copies of this basic construction to build the figure $\nabla$ as the image of a local isometry from a subset of $\euc^2$. We'll show that by patching appropriately, this subset can be a rectangle with the necessary edge identifications from Figure \ref{fig:klein-bottle-rectangle}.

 Indeed, by setting $\theta=\frac{\pi}{3}$ and restricting the constructed map $f$ to the rectangle $D_0$ where $-\pi r\leq u\leq \pi r$ and $-\frac{s}{2}\leq v\leq \frac{s}{2}$, the image of $f$ is one-third of the figure $\nabla$. Changing coordinates appropriately, one can construct a map $\phi$ from the rectangle consisting of three stacked copies of $D_0$ to three copies of the image of $f$ glued together to form $\nabla$.
 
   \begin{figure}[h!tbp]
\centering
\begin{minipage}[b]{0.45\textwidth} 
 \begin{center}
	\begin{tikzpicture}[scale=0.6]
	 \def\a{0.577}; % 1/sqrt(3)
   \def\b{1.732}; % sqrt(3)
   \fill[green!25!white] (-pi,-2) rectangle (pi,2);
	 \draw[dashed] (-pi,0) sin (-pi/2,-\a) cos (0,0) sin (pi/2,\a) cos (pi,0);
	 
	 \begin{scope}[yshift=4cm,xscale=-1]
	  \fill[blue!25!white] (-pi,-2) rectangle (pi,2);
	 \draw[dashed] (-pi,0) sin (-pi/2,-\a) cos (0,0) sin (pi/2,\a) cos (pi,0);
	 \end{scope}
	 
	 \begin{scope}[yshift=-4cm,xscale=-1]
	  \fill[red!25!white] (-pi,-2) rectangle (pi,2);
	 \draw[dashed] (-pi,0) sin (-pi/2,-\a) cos (0,0) sin (pi/2,\a) cos (pi,0);
	 \end{scope}
	 
	 \draw[thick,->,cyan] (pi,2) -- (0,2);
	 \draw[thick,cyan] (0,2) -- node[below=0.5cm,green]{$D_0$} (-pi,2);
	 
	 \draw[thick,->,magenta] (pi,-6) -- (0,-6);
	 \draw[thick,magenta] (0,-6) -- (-pi,-6);
	 
	 \draw[thick,->,magenta] (-pi,6) -- (0,6);
	 \draw[thick,magenta] (0,6) -- node[below=0.5cm,blue]{$D_1$} (pi,6);
	 
	 \draw[thick,->,yellow] (-pi,-2) -- (0,-2);
	 \draw[thick,yellow] (0,-2) -- node[below=0.5cm,red]{$D_{-1}$} (pi,-2);
	 
	 \draw[thick, ->] (-pi,-6) -- (-pi,-0.3);
	 \draw[thick, ->] (-pi,-0.3) -- (-pi,-0.2);
	 \draw[thick] (-pi,6) -- (-pi,-0.2);
	 
	 \draw[thick, ->] (pi,-6) -- (pi,-0.3);
	 \draw[thick, ->] (pi,-0.3) -- (pi,-0.2);
	 \draw[thick] (pi,6) -- (pi,-0.2);
	 
	\end{tikzpicture}
	
        	\subcaption{} 
        	\label{fig:phi-domain}
	\end{center}
\end{minipage}
\begin{minipage}[b]{0.45\textwidth} 
 \begin{center}
	\includegraphics[height=2.5in]{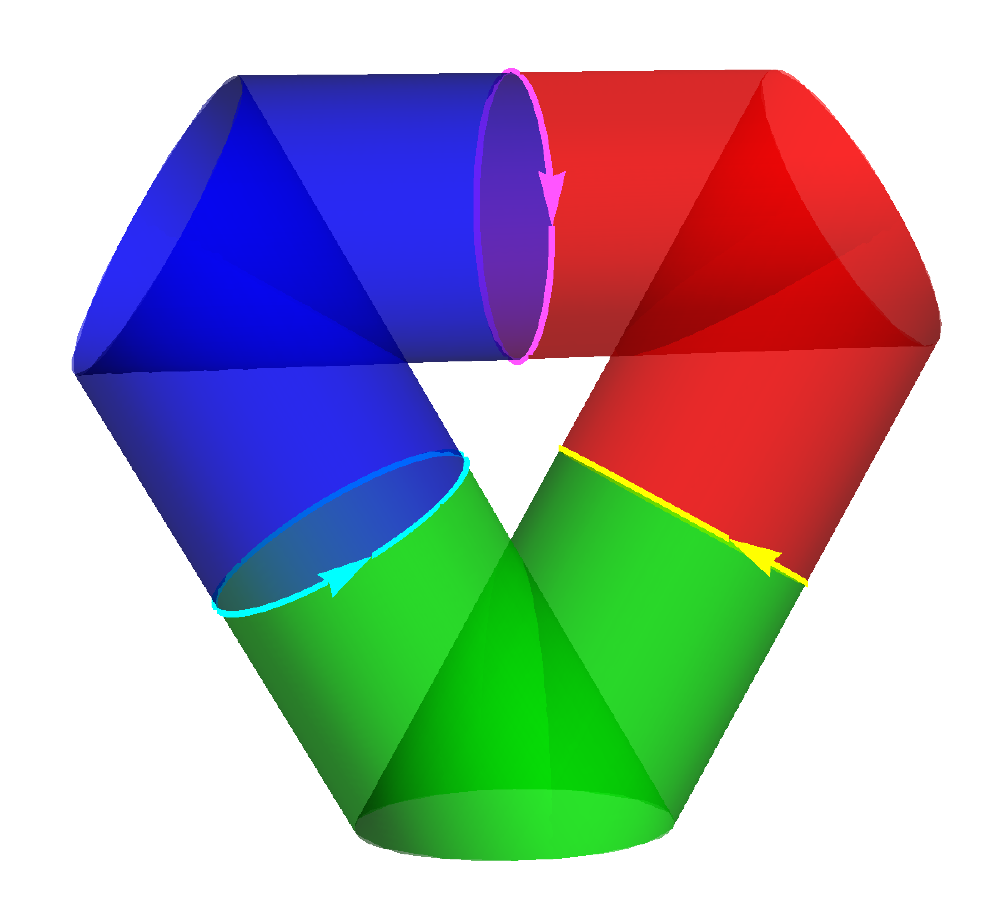}
        	\subcaption{}
        	\label{fig:glued-nabla}
 \end{center}
\end{minipage}

\caption{Gluing copies of the basic construction. The rectangular domain (a) consists of stacked copies of $D_0$; the creases are shown as dotted lines and the oriented gluing edges are marked with arrows. Glued appropriately, we obtain the figure $\nabla$ in (b).}
\label{fig:gluing}
\end{figure}
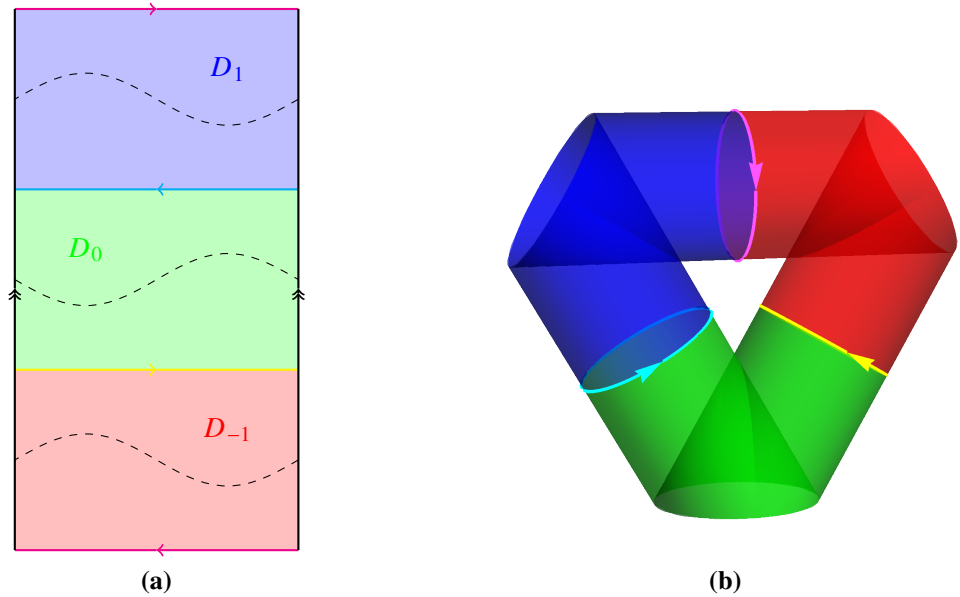
 
 We note that it is clear from the formula (\ref{eqn:cylinder-formula}) that $f_1$ respects the identification of the vertical edges, and therefore $\phi$ does as well. By construction, the horizontal edges both get isometrically mapped to the same circle, and since the handedness of the cylinder parameterization is switched an odd number of times, the orientation of these maps must be opposite one another, as required. In this way, $\phi$ only fails to be injective along the three semi-elliptical pairwise intersections of the cylinders highlighted in figure \ref{fig:joint}. The rest of the claimed properties of $\phi$ from the Introduction follow from the basic construction: appropriate gluing ensures that $\phi$ is everywhere continuous, $\phi$ is smooth except along the three copies of $\alpha$ (which map to the folded ellipses), and $\phi$ is a local isometry except at each of the two inflection points of each  copy of $\alpha$.

We note that similarly taking the union of appropriate pieces cylinders of equal radii centered around the edges of any $n$-sided polygon in $\mathcal E^3$ will result in a figure that is the image of torus if $n$ is even or a Klein bottle if $n$ is odd.

  \begin{figure}[h!tbp]
\centering
\begin{minipage}[b]{0.45\textwidth} 
 \begin{center}
	\includegraphics[height=2in]{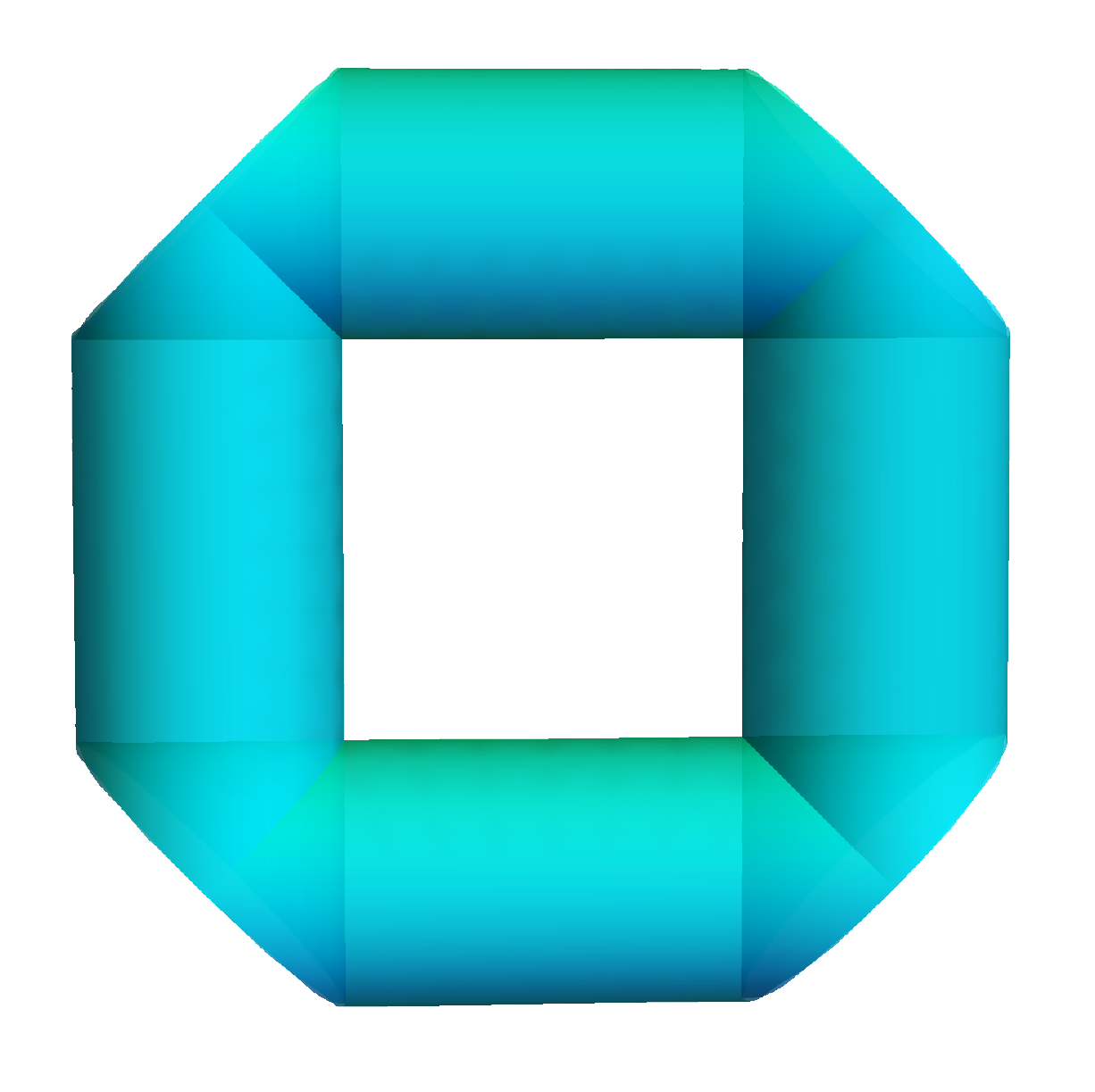}
        	\subcaption{} 
        	\label{fig:nabla-4}
	\end{center}
\end{minipage}
\begin{minipage}[b]{0.45\textwidth} 
 \begin{center}
	\includegraphics[height=2in]{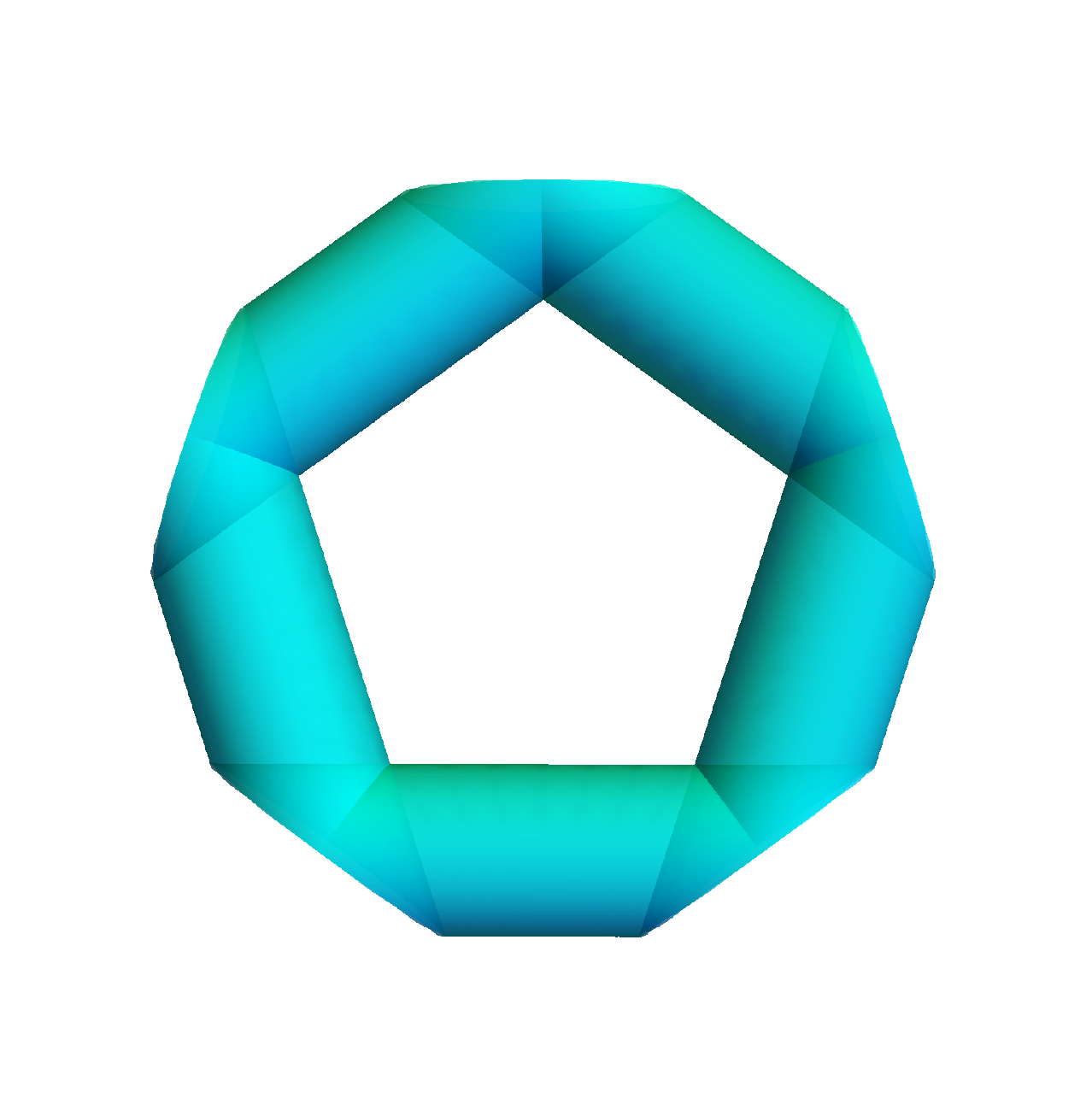}
        	\subcaption{}
        	\label{fig:nabla-5}
 \end{center}
\end{minipage}

\caption{By appropriately gluing $n$ copies of the basic construction, one obtains (a) an image of the torus when $n$ is even or (b) a Klein bottle when $n$ is odd.}
\label{fig:nablas}
\end{figure}

\section*{Illustration Through Sculpture}\label{section:medium}

The aim of the present sculpture is largely a pedagogical one: to make concrete the consequences of concepts taught in courses on topology and geometry. The Klein bottle is a familiar object to mathematicians in the abstract, and while sculptures and drawings like the one in Figure \ref{fig:usual-klein-bottle} help one understand its topology, they do not attempt to illustrate its admission of flat geometry.

I made the sculpture from clear, printed transparency film by cutting, creasing, and taping. While the creases are meant to be seen, the cuts and tape are hidden, giving the illusion of physical continuity to illustrate the underlying mathematical continuity. The grid pattern and transparent medium are further meant to aide the viewer in willfully looking past the self-intersections and thereby focus on the properties of the sculpture inherited from the abstract Klein bottle. In particular, the cylindrical pieces are meant to appear to pass through each other---their semi-elliptical intersection (achieved by cutting one of the cylinders to allow the other to pass through) easily ignored.

\begin{figure}[h!tbp]
\centering
\begin{minipage}[b]{0.55\textwidth} 
 \begin{center}
	\includegraphics[height=2in]{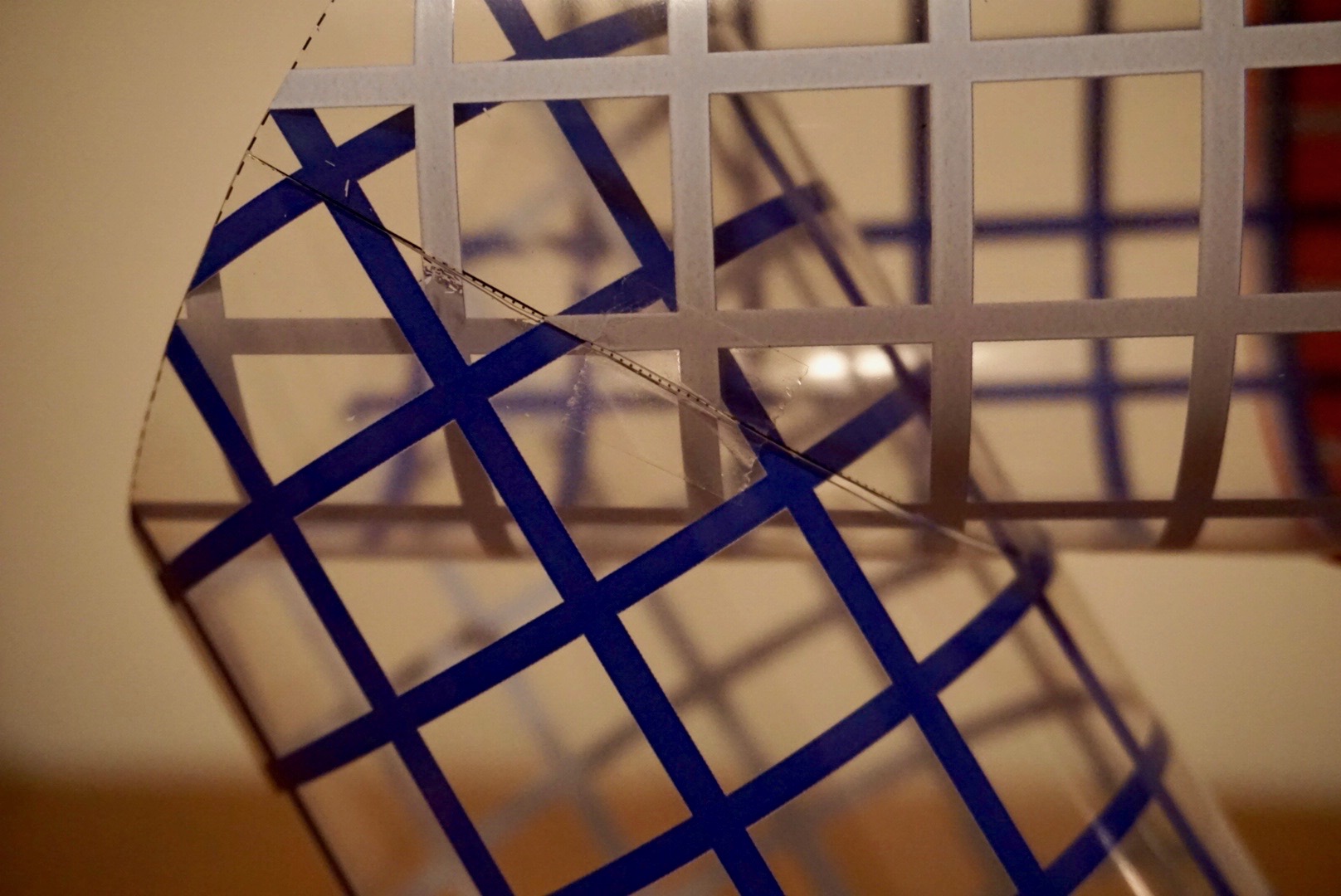}
        	\subcaption{} % Add subcaption text if desired, or use \subcaption* to suppress (a), (b), etc. labels
        	\label{fig:pass-through}
	\end{center}
\end{minipage}
\begin{minipage}[b]{0.4\textwidth} 
 \begin{center}
	\includegraphics[height=2in]{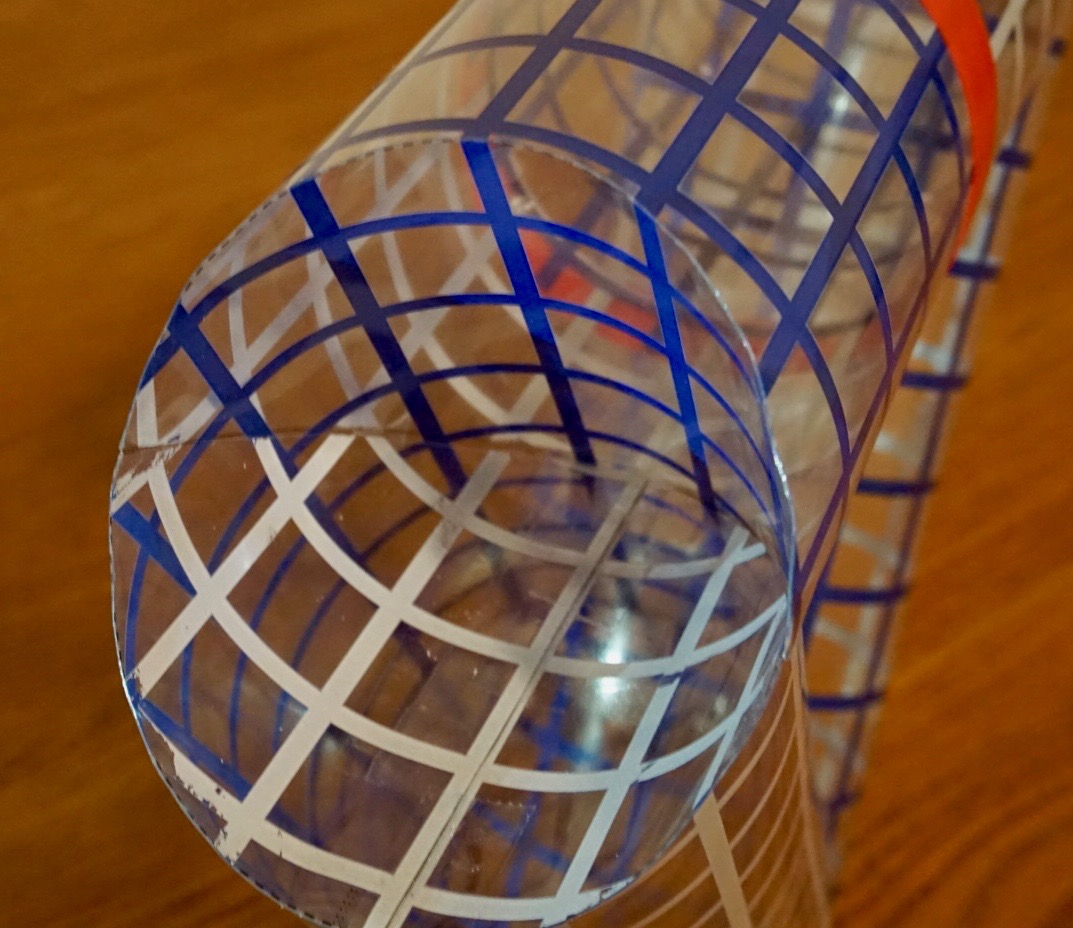}
        	\subcaption{} % Add subcaption text if desired, or use \subcaption* to suppress (a), (b), etc. labels
        	\label{fig:self-intersection}
 \end{center}
\end{minipage}

\caption{Details of a self-intersection (a) and folded ridge (b) on the sculpture.}
\label{fig:sculpture-details}
\end{figure}

%\begin{figure}[H]
% \begin{center}
%  \includegraphics[width=3in]{images/pass-through}
% \end{center}
%\end{figure}

Consideration of these self-intersections suspended, the sculpture apparently represents a compact genus-one non-orientable surface (the non-orientability emphasized by the coloring of the grid pattern). A viewer with some background in topology will quickly identify this as the Klein bottle, even if this particular manifestation is unfamiliar.

A viewer may then recognize that the construction from inelastic transparency film demonstrates the flat geometry of the surface. As already noted, compact surfaces with topology different from the Klein bottle and torus do not admit flat geometry. The use of origami, including curved-crease origami, to model surfaces with flat metrics has been studied before \cite{huffman76,ft99}, and in this way, the medium itself illustrates the flat metric on the Klein bottle. The grid pattern is designed to further emphasize this flatness. As noted earlier, a ``texture'' on a flat plane can be mapped undistorted onto a flat surface, even across folds and vertices, as it is in the present sculpture.

%\begin{figure}[H]
% \begin{center}
%  \includegraphics[width=3in]{images/ridge}
% \end{center}
%\end{figure}

\section*{Mathematical and Artistic Context}

The context of the present sculpture of the flat Klein bottle and its mathematical underpinnings is best understood by exploring the more well-developed set of works regarding the flat torus.

\begin{figure}[h!tbp]
\centering
\begin{minipage}[b]{0.3\textwidth} 
 \begin{center}
	\includegraphics[height=1.1in]{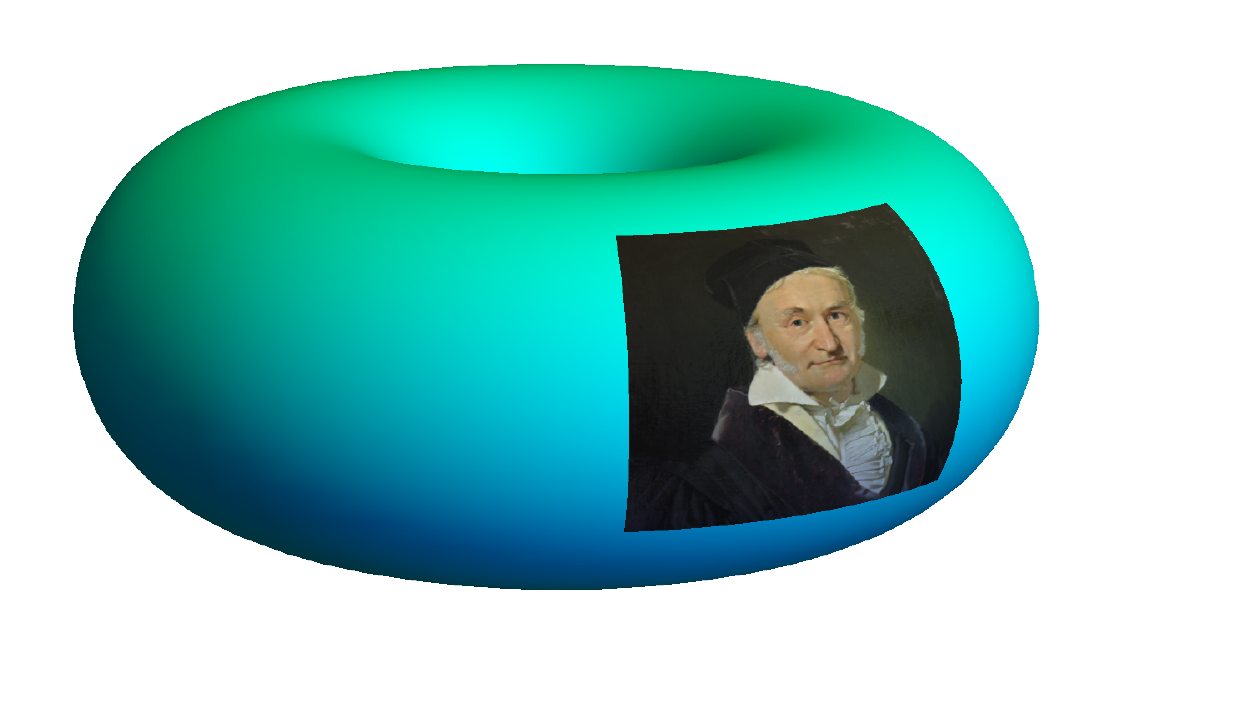}
        	\subcaption{} 
        	\label{fig:gauss-torus}
	\end{center}
\end{minipage}
\begin{minipage}[b]{0.3\textwidth} 
 \begin{center}
	\includegraphics[height=1.5in]{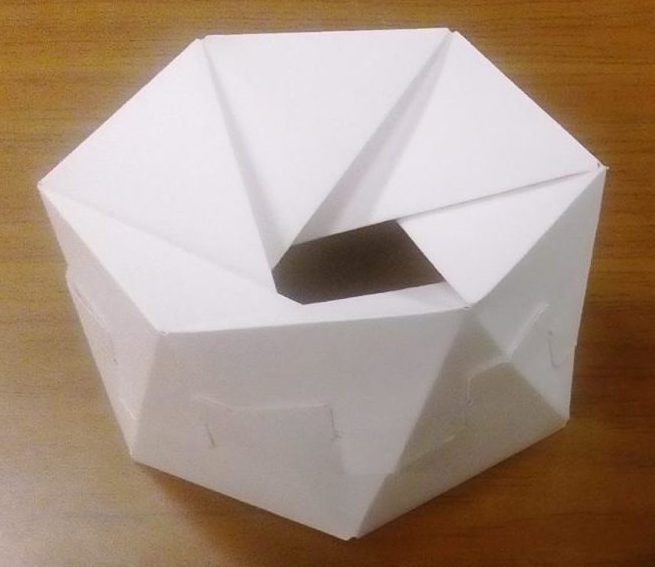}
        	\subcaption{}
        	\label{fig:paper-torus}
 \end{center}
\end{minipage}
\begin{minipage}[b]{0.3\textwidth} 
 \begin{center}
	\includegraphics[height=1.5in]{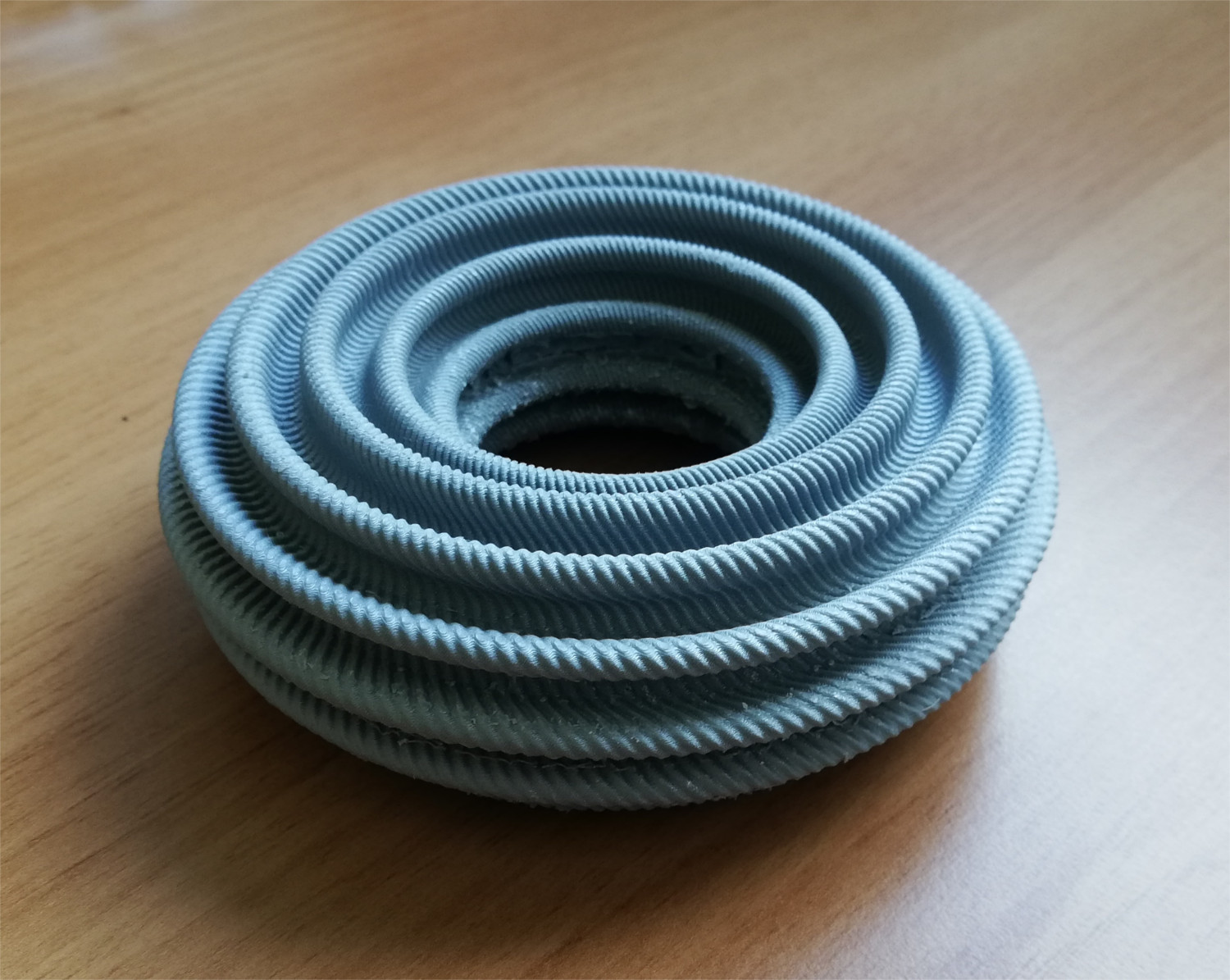}
        	\subcaption{} % Add subcaption text if desired, or use \subcaption* to suppress (a), (b), etc. labels
        	\label{fig:c1-torus}
 \end{center}
\end{minipage}

\caption{Embeddings of the flat torus: (a) A smooth embedding of the torus cannot be flat. (b) An origami sculpture of a piecewise-linear flat torus \cite{paperflattori}. (c) A 3D printed sculpture of a continuously differentiable flat torus \cite{tm13}.}
\label{fig:tori}
\end{figure}

As we mentioned earlier, any flat compact surface necessarily has the global topology of either a Klein bottle or torus. Unlike the Klein bottle, the torus may be smoothly embedded in euclidean three-space. However, no smooth embedding of the torus can also preserve the geometry of a flat torus; that is, it cannot be both smoothly and isometrically embedded (notice the necessary distortion in Figure \ref{fig:gauss-torus}).
 
 However, explicit examples of continuous piecewise-linear isometries of the flat torus into euclidean three-space do exist, as illustrated by Segerman using 3D printing in \cite{segerman16} and M\'alaga and Leli\`evre in \cite{paperflattori} using origami.
 
 In fact, Nash \cite{nash54} and Kuiper \cite{kuiper55} proved theorems in the 1950s which show that there exists a continuously differentiable (but not twice-differentiable) isometric embedding of the flat torus into euclidean three-space, although its image is necessarily a fractal. An explicit example of such an embedding was not found until 2012 by Borrelli \emph{et al} in \cite{bjlt12}, and a 3D printed model, pictured in Figure \ref{fig:c1-torus}, of the embedded flat torus was presented by Henocque and Marin in \cite{tm13}.
 
\begin{figure}[h!tbp]
\centering
\begin{minipage}[b]{0.45\textwidth} 
 \begin{center}
	\includegraphics[height=2.5in]{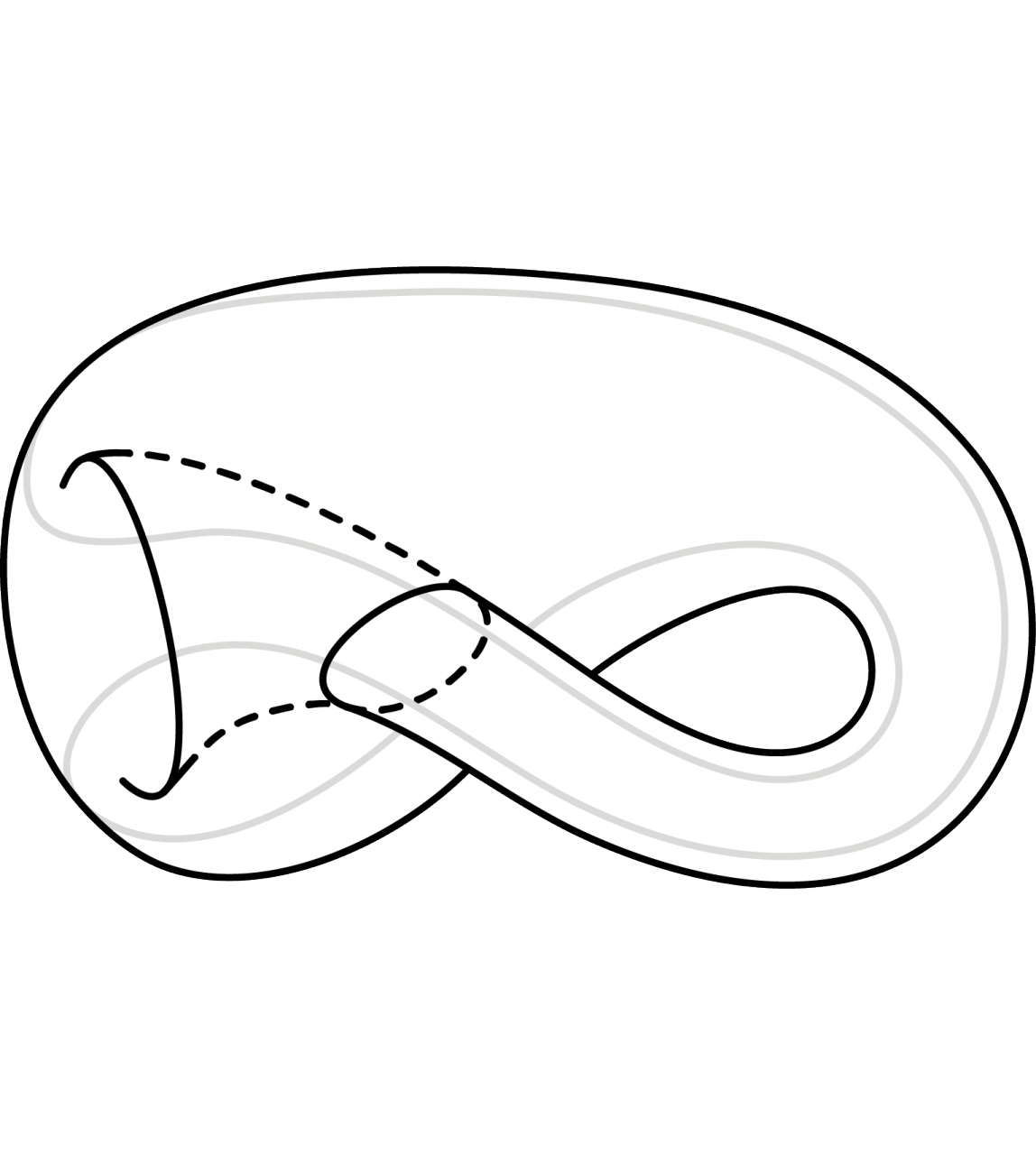}
        	\subcaption{} 
        	\label{fig:usual-klein-bottle}
	\end{center}
\end{minipage}
\begin{minipage}[b]{0.45\textwidth} 
 \begin{center}
	\includegraphics[height=2.5in]{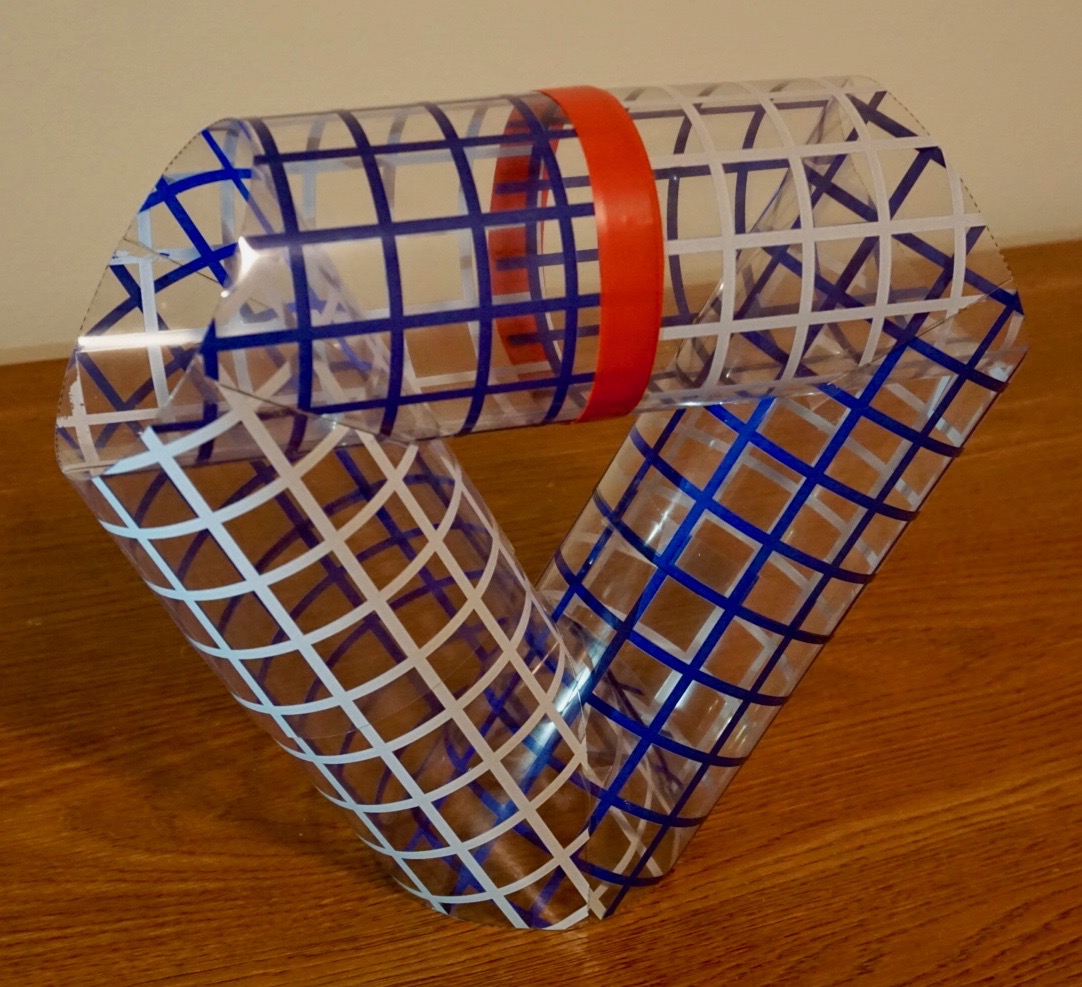}
        	\subcaption{}
        	\label{fig:nabla-color}
 \end{center}
\end{minipage}

\caption{Images of the Klein bottle: (a) A smooth local embedding of the Klein bottle cannot be flat. (b) A sculpture of $\nabla$.}
\label{fig:klein}
\end{figure}

 One can view the present sculpture as filling in a piece of a parallel and mostly untold story about the Klein bottle. The Klein bottle cannot be globally embedded in euclidean three-space, but it can be smoothly locally embedded as in Figure \ref{fig:usual-klein-bottle}. The results of Nash and Kuiper show that smoothness can be traded for local isometry and continuous differentiability, but no explicit example of such an isometry has been exhibited. We can think of the map $\phi$ constructed here as satisfying a weaker of the set of properties than those guaranteed by Nash and Kuiper. Thus the construction of $\nabla$ and its sculpture are presented as a first step towards developing illustrations of the flat Klein bottle that parallel those exhibited for the flat torus.
 
%%%%%%%%%%%%%
\section*{Summary and Conclusions}

While the existence of a flat Klein bottle in the abstract is well-known to mathematicians, I hope that by exhibiting its image in euclidean three-space with much of the topology and geometry still intact will help make this fact feel more concrete, especially to those encountering such ideas for the first time.

%%%%%%%%%%%%%%%%%%%%%%%%%%%%%%%%%%%%%%%
\section*{Acknowledgements}

This material is based upon work supported by the National Science Foundation under Grant No. DMS-1439786 and the Alfred P. Sloan Foundation award G-2019-11406 while the author was in residence at the Institute for Computational and Experimental Research in Mathematics in Providence, RI, during the Illustrating Mathematics program. Along the many other artists and mathematicians I met and worked with at ICERM whose ideas worked their way into this project, I wish to thank John Edmark for introducing me to curved-crease origami, Glen Whitney introducing me to the origami flat torus.

%%%%%%%%%%%%%%%%%%%%%%%%%%%%%%%%%%%%%%%
% References %
    
{\setlength{\baselineskip}{13pt} % tighten line spacing for bibliography
\raggedright				% no right justification for References

\end{document}